\documentclass{article}
\usepackage{amsfonts, mathrsfs, amsmath, amssymb}
\newtheorem{def1}{Definition}[section]
\newtheorem{pro}{Proposition}[section]
\newtheorem{rem}{Remark}[section]

\newtheorem{cor} {Corollary}[section]
\begin{document}
\newtheorem{thm}{Theorem}[section]
\title{The categories $L\textbf{-Top}_{0}$ and $L\textbf{-Sob}$ as epireflective hulls}
\author {Rana Noor\thanks{rananoor4@gmail.com}\\ {\em Department of Mathematics},\\ {\em Banaras Hindu University},\\ {\em Varanasi-221 005, India}\\\\Arun K. Srivastava\thanks{arunksrivastava@gmail.com} \\ {\em Department of Mathematics}\\ {\em and}\\ {\em Centre for Interdisciplinary Mathematical Sciences},\\ {\em Banaras Hindu University},\\ {\em Varanasi-221 005, India}} 
\date{}
\maketitle
\abstract
We show that the category $L\textbf{-Top}_{0}$ of $T_{0}$-$L$-topological spaces is the epireflective hull of Sierpinski $L$-topological space in the category $L\textbf{-Top}$ of $L$-topological spaces and the category $L\textbf{-Sob}$ of sober $L$-topological spaces is the epireflective hull of Sierpinski $L$-topological space in the category $L\textbf{-Top}_{0}$.
\section{Introduction}
\qquad Rodabaugh \cite{Roda}  extended the notion of sobriety in topology (cf. \cite{John}) to $L$-topology. It is known that the category of sober topological spaces is the epireflective hull of the two-point Sierpinski topological space in the category $\textbf{Top}_{0}$ of $T_{0}$-topological spaces (cf. \cite{Nel}). H\"{o}hle \cite{Hole} has already pointed out that the category of sober $L$-topological spaces is reflective in the category $L\textbf{-Top}$ of $L$-topological spaces. We show here that the category of sober $L$-topological spaces is in fact the epireflective hull of the `Sierpinski $L$-topological space' in the category $L\textbf{-Top}_{0}$ of $T_{0}$-$L$-topological space.
\section{Preliminaries}
\qquad For the category theoretic notions used here, \cite{AHS} may be referred. All subcategories are assumed to be full and replete.\\

Throughout this paper, $L$ denotes a frame (a complete lattice which satisfies the first infinite distributive law: $a\wedge(\vee b_{i})= \vee(a\wedge b_{i})$), with $0$ and $1$ being its least and largest elements respectively. A \textit{frame map} is a lattice homomorphism between frames which preserves finite infima and arbitrary suprema. A subset $A\subseteq L$ of a frame $L$ is called a \textit{subframe} if $A$ is a frame under the order induced by $L$.\\
We recall some definitions to make the paper self contained.




For a given set $X$,  $L^{X}$ is also a frame under the order induced by $L$, having $\bar{0}$ and $\bar{1}$ as its least and largest elements respectively. Members of $L^{X}$, are known as $L$-sets in $X$ \cite{Goguen}.
\begin{def1}{\rm \cite{Hole}}
A family $\tau$ of $L$-sets in a set $X$ is called an {\rm $L$-topology} on $X$, and the pair $(X,\tau)$ an {\rm $L$-topological space}, if $\tau$ is a subframe of $L^{X}$. Furthermore, a map $f:(X,\tau)\rightarrow (Y,\delta)$ between two $L$-topological spaces is called {\rm continuous} if $f^{\leftarrow}(\nu)\in\tau$, for every $\nu\in\delta$.  
\end{def1}

Let $L\textbf{-Top}$ denote the category of all $L$-topological spaces and their continuous maps.

Given an $L$-topological space $(X,\tau)$ and a subset $Y\subseteq X$, $\tau_{Y}=\{\mu\wedge 1_{Y}$ $|$ $\mu\in\tau\}$ is an $L$-topology on $Y$. The $L$-topological space $(Y,\tau_{Y})$ is called a \textit{subspace} of $(X,\tau)$.

The notions of 
\textit{open map}, \textit{homeomorphism} and \textit{embedding} in $L\textbf{-Top}$ are on expected lines.
\begin{def1}{\rm \cite{AHS}}
A concrete category $(\mathscr C,U)$ over a category $\mathscr A$ is said to be {\rm uniquely transportable} if for every $X\in ob\mathscr C$ and every $\mathscr A$-isomorphism $k:U(X)\rightarrow A$ there exists a unique $Y\in ob\mathscr C$ with $U(Y)=A$ such that $k':X\rightarrow Y$ is a $\mathscr C$-isomorphism with $U(k')=k$.
\end{def1}

\begin{def1}{\rm \cite{Mane}}
A {\rm category $\mathscr C$ of sets with structures} is a uniquely transportable construct $(\mathscr C,U)$ in which
a class $C_{\mathscr C}(X)$ of `$\mathscr C$-structures' is assigned to every set $X$, which is in one-to-one correspondence with the ‘fibre’ $\{A\in ob\mathscr C$ $|$ $U(A) = X\}$ of $X$; thus each $\mathscr C$-object $A$ can be identified with a pair $(X, s)$ $($and called a ‘set with a structure’$)$, where $U(A) = X$ and $s\in C_{\mathscr C}(X)$.
\end{def1}
\begin{def1}{\rm \cite{AHS}}
Given a category $\mathscr C=(\mathscr C,U)$ of sets with structures,
\begin{itemize}
\item a family ${\mathscr F}=\{f_i:(X, s)\rightarrow (Y_i, t_i)$ $|$ $i\in I\}$ of $\mathscr C$-morphisms is said to be {\rm initial} $($or {\rm optimal}, as in {\rm \cite{Mane}}$)$ if for every $\mathscr C$-structured set $(Z, u)$ and a map $g:Z\rightarrow X$,  $g:(Z, u)\rightarrow (X, s)$ if $f_i\circ g:(Z, u)\rightarrow (X, s),$ for every $ i\in I$.
\item if for a family ${\mathscr F}=\{f_i:X\rightarrow U(Y_i, t_i)$ $|$ $i\in I\}$ of maps, where $X$ is a set and each $(Y_i, t_i)$ is a $\mathscr C$-structured set, there exists $s\in C_{\mathscr C}(X)$ such that the family $\{f_i:(X, s)\rightarrow (Y_i, t_i)$ $|$ $i\in I\}$ is optimal, then $s$ is called an {\rm initial lift} $($or {\rm optimal lift}, as in {\rm \cite{Mane}}$)$ of the family ${\mathscr F}$.
\end{itemize}
\end{def1}
\begin{def1}{\rm \cite{Mane}}
Given a category $\mathscr C=(\mathscr C,U)$ of sets with structures, an object $S$ of $\mathscr C$ is called a {\rm Sierpinski object} if for every $X=(X,s)\in ob\mathscr C$, the family of all $\mathscr C$-morphisms from $X$ to $S$ is initial.
\end{def1}

The category $L\textbf{-Top}$ is obviously a category of sets with structure. Moreover it is a topological construct \cite{Hole}. We shall use $V$ to denote the forgetful functor from $L\textbf{-Top}$ to the category of sets and functions, wherever needed.\\
\begin{rem}
Let $\zeta$ be a family of $L$-sets in a set $X$. Let $\xi$ be the collection of all finite infima of members of $\zeta$ and $\tau$ be the collection of all suprema of members of $\xi$. It can be verified that $\tau$ is an $L$-topology on $X$. We shall denote it as $<\zeta>$.
\end{rem}
\begin{def1}
Let $\mathscr F=\{f_{i}:X\rightarrow (Y_{i},\delta_{i})$ $|$ $i\in I\}$ be a family of maps, where $X$ is a set and $\{(Y_{i},\delta_{i})$ $|$ $i\in I\}$ is a family of $L$-topological spaces. Then $<\{f_{i}^{\leftarrow}(\mu)$ $|$ $\mu\in\delta_{i}, i\in I\}>$ is called the {\rm initial $L$-topology} on $X$ induced by $\mathscr F$.
\end{def1}
\begin{def1}
Given a family $\{(X_{i},\tau_{i})$ $|$ $i\in I\}$ of $L$-topological spaces, the initial $L$-topology on $X$ $(=\displaystyle\prod_{i\in I}X_{i})$ induced by all projection maps $\{p_{i}:X\rightarrow (X_{i},\tau_{i})$ $|$ $i\in I\}$ is called the {\rm product $L$-topology}.
\end{def1}

Let $(X,\tau)$ be an $L$-topological space, $Y$ be a set and $f:X\rightarrow Y$ be a surjective map. Then $\tau/f =\{\nu\in L^{Y}$ $|$ $f^{\leftarrow}(\nu)\in\tau\}$ is an $L$-topology on $Y$, called the \textit{quotient $L$-topology}, and the pair $(Y,\tau/f)$ is called the \textit{quotient $L$-topological space} with respect to $(X, \tau)$ and $f$. The resulting continuous map $f:(X,\tau)\rightarrow (Y,\tau/f)$ is called a \textit{quotient map}.\\

Let $\mathscr C$ be a category and let $\mathscr H$ and epi$\mathscr{C}$, respectively, denote some class of $\mathscr C$-morphisms and the class of all $\mathscr C$-epimorphisms.
\begin{def1}
A $\mathscr C$-object $X$ is called $\mathscr H$-{\rm injective} if for every $e:Y\rightarrow Z$ in $\mathscr H$ and every $\mathscr C$-morphism $f:Y\rightarrow X$, there exists a $\mathscr C$-morphism $g:Z\rightarrow X$ such that $g\circ e=f$.
\end{def1}

\begin{def1}
A $\mathscr C$-object $X$ is called a {\rm cogenerator} $($called a {\rm coseparator} in {\rm \cite{AHS}}$)$ in $\mathscr C$ if for every pair of distinct $f,g\in \mathscr C (Y,Z)$, there exists $h\in \mathscr C (Z,X)$ such that $h\circ f\neq h\circ g$.
\end{def1}
Let $\mathscr R$ be a subcategory of $\mathscr C$.
\begin{def1}{\rm \cite{Giu}}
$\mathscr R$ is said to be {\rm epireflective} in $\mathscr C$ if for each $\mathscr C$-object $X$, there exists an epimorphism $r_{X}:X\rightarrow RX$, with $RX\in ob\mathscr {R}$, such that for each $\mathscr C$-morphism $f:X\rightarrow Y$, with $Y\in ob\mathscr {R}$, there exists a unique $\mathscr {R}$-morphism $f^{*}:RX\rightarrow Y$, such that $f^{*}\circ r_{X}=f$. If moreover, each $r_{X}\in \mathscr {H}$ and $f^{*}$ is a $\mathscr C$-isomorphism, whenever $f\in epi\mathscr{C}\bigcap \mathscr H$, then $\mathscr R$ is said to be an {\rm $\mathscr {H}$-firm epireflective subcategory} of $\mathscr C$.
\end{def1}
\section{The Sierpinski $L$-topological Space}
\qquad Consider the frame $L$. Then $<id_{L}>$, where $id_{L}$ is the identity map on $L$, is an $L$-topology on $L$. Call $(L, <id_{L}>)$ the \textit{Sierpinski $L$-topological space} and denote it as $L_{S}$.\\\\
The following result is easy to verify.
\begin{pro}
Let $(X,\tau)\in ob L\textbf{-Top}$. Then $\mu\in\tau$ iff $\mu:(X,\tau)\rightarrow L_{S}$ is continuous.
\end{pro}
\begin{thm}
$L_{S}$ is a Sierpinski object in $L\textbf{-Top}$.
\end{thm}
\textbf{Proof}: Let $(X,\tau)\in obL\textbf{-Top}$ and consider the family $\mathscr F=\{f:(X,\tau)\rightarrow L_{S}$ $|$ $f$ is continuous$\}$. Let $(Y,\delta)\in ob L\textbf{-Top}$ and let $g:Y\rightarrow X$ be a map such that $f\circ g:(Y,\delta)\rightarrow L_{S}$ is continuous, for every  $f\in\mathscr F$. We have to show that $g:(Y,\delta)\rightarrow (X,\tau)$ is continuous. Let $\mu\in\tau$. Then by Proposition $3.1$, $\mu:(X,\tau)\rightarrow L_{S}$ is continuous and hence $\mu\in\mathscr F$. So, $\mu\circ g:(Y,\delta)\rightarrow L_{S}$ is continuous. Again by Proposition $3.1$, $\mu\circ g\in\delta$. But $g^{\leftarrow}(\mu)=\mu\circ g$, implying that $g$ is continuous. Thus $L_{S}$ is a Sierpinski object in $L\textbf{-Top}$. $\Box$
\begin{def1}{\rm \cite{Roda}}
An $L$-topological space $(X,\tau)$ is called $T_{0}$ if for every distinct $x,y\in X$, there exists some $\mu\in\tau$ such that $\mu(x)\neq\mu(y)$.
\end{def1}
The Sierpinski $L$-topological space $L_{S}$ is $T_{0}$.\\

Let $L$\textbf{-Top}$_{0}$ denote the subcategory of $L$\textbf{-Top}, whose objects are $T_{0}$-$L$-topological spaces.\\


From now on, we write `injective' in place of `$\mathscr H$-injective', when $\mathscr H$ is the class of all embeddings in $L\textbf{-Top}_{0}$.
\begin{thm}
$L_{S}$ is an injective cogenerator in $L\textbf{-Top}_{0}$.
\end{thm}
\textbf{Proof}: Let $e:(X,\tau)\rightarrow (Y,\delta)$ be an embedding in $L\textbf{-Top}_{0}$ and let $f:(X,\tau)\rightarrow L_{S}$ be a continuous map. Then $f\in\tau$, implying that $e^{\rightarrow}(f)\in \delta_{e(X)}$. So, there exists $\nu\in\delta$ such that $e^{\rightarrow}(f)=\nu\wedge 1_{e(X)}$ and $\nu\circ e=f$. As $\nu\in\delta$, $\nu:(Y,\delta)\rightarrow L_{S}$ is continuous. Thus $L_{S}$ is injective.

Consider distinct pair of morphism $f,g:(X,\tau)\rightarrow (Y,\delta)$ in $L\textbf{-Top}_{0}$. Then there exists $x\in X$ such that $f(x)\neq g(x)$. As $(Y,\delta)$ is $T_{0}$, there exists $\nu\in\delta$ such that $\nu(f(x))\neq\nu(g(x))$ implying that $\nu\circ f\neq \nu\circ g$. As $\nu\in\delta$, $\nu:(Y,\delta)\rightarrow L_{S}$ is continuous. Thus $L_{S}$ is a cogenerator. $\Box$
\section{A characterization of $L$-Top}
\qquad Manes \cite{Mane} obtained a characterization of the category $\textbf{Top}$ of topological spaces with the help of usual two-point Sierpinski topological space $2_{S}$. Srivastava \cite{Arun} obtained a characterization of the category $\textbf{FTop}$ of fuzzy topological spaces with the help of fuzzy Sierpinski space $I_{S}$. In this section, we shall give a characterization of $L\textbf{-Top}$ with the help of $L_{S}$.
\begin{thm}
A category $(\mathscr C,U)$ of sets with structures is concretely isomorphic to $L\textbf{-Top}$ iff there exists $(S,u)\in ob\mathscr C$ with the underlying set $S=L$, satisfying the following condition:
\begin{enumerate}
\item $(S,u)$ is a Sierpinski object in $\mathscr C$,
\item every family $\{f_{i}:X\rightarrow U(S,u)$ $|$ $i\in I\}$ has an initial lift,
\item the map $sup:(S,u)^{X}\rightarrow (S,u)$ is $\mathscr C$-morphism for every set $X$ $($here $(S,u)^{X}=(S^{X},t)$, where $t$ is the initial lift of all projection from $S^{X}$ to $U(S,u))$,
\item The map $inf: (S,u)^{X}\rightarrow (S,u)$ is $\mathscr C$-morphism for every finite set $X$,
\item for every $(X,s)\in ob\mathscr C$ and for every initial family $\mathscr F$ of $\mathscr C$-morphisms from $(X,s)$ to $(S,u)$, the initial lift of $\{f:X\rightarrow U(S,u)$ $|$ $f\in\mathscr F\}$ contains every $\mathscr C$-morphism $g:(X,s)\rightarrow (S,u)$. 
\end{enumerate}
\end{thm}
\textbf{Proof}: Let $(\mathscr C,U)$ be concretely isomorphic to $L\textbf{-Top}$. We show that $L\textbf{-Top}$ satisfies $(1)-(5)$. In view of Theorem $3.1$, $L_{S}$ is a Sierpinski object in $L\textbf{-Top}$. So, $L\textbf{-Top}$ satisfies $(1)$. It is easy to see that every family $\mathscr G=\{f_{i}:X\rightarrow V(L_{S})$ $|$ $i\in I\}$ in $L\textbf{-Top}$ has an initial lift, being the initial $L$-topology on $X$ induced by $\mathscr G$. So, $L\textbf{-Top}$ satisfies $(2)$. Let $\mu\in L_{S}^{X}$ and for $x\in X$, $p_{x}$ denote the $x$-th projection map from $L^{X}$ to $V(L_{S})$. Then $sup (\mu)= \displaystyle\bigvee_{x\in X}\mu(x)=\displaystyle\bigvee_{x\in X}p_{x}(\mu)$ implying that $sup=\vee\{p_{x}$ $|$ $x\in X\}$. As $p_{x}$ is open in $L_{S}^{X}$, $sup$ is also open in $L_{S}^{X}$. Thus $sup:L_{S}^{X}\rightarrow L_{S}$ is continuous, for every set $X$. Similarly, $inf=\wedge\{p_{x}$ $|$ $x\in X\}$ and therefore $inf$ is also continuous, for every finite set $X$. Hence $L\textbf{-Top}$ satisfies $(3)$ and $(4)$. Let $(X,\tau)\in ob L\textbf{-Top}$ and $\mathscr F=\{f$ $|$ $f:(X,\tau)\rightarrow L_{S}$ is continuous$\}$ be an initial family. Then $\mathscr F\subseteq \tau$. Let $\sigma$ be the initial lift of $\{f:X\rightarrow V(L_{S})$ $|$ $f\in\mathscr F\}$. Then $\sigma\subseteq\tau$ and $f:(X,\sigma)\rightarrow L_{S}$ is continuous, for every $f\in \mathscr F$. Consider the identity map $id:(X,\sigma)\rightarrow (X,\tau)$. As $f\circ id=f$, for every $f\in\mathscr F$, $id$ is continuous and hence $\tau\subseteq \sigma$. Thus $L\textbf{-Top}$ satisfies $(5)$.

Conversely, let $(\mathscr C,U)$ satisfies the given conditions. We show that $(\mathscr C,U)$ and $L\textbf{-Top}$ are concretely isomorphic. For showing this we have to produce two concrete functors $F:\mathscr C\rightarrow L\textbf{-Top}$ and $G:L\textbf{-Top}\rightarrow\mathscr C$ which are inverses to each other.\\ 
Let $(X,s)\in ob\mathscr C$ and let $\tau_{s}$ be the initial lift of the family $\mathscr F=\{f:X\rightarrow V(L_{S})$ $|$ $f:(X,s)\rightarrow (S,u)$ is $\mathscr C$-morphism$\}$. Then $\tau_{s}=<\{f^{\leftarrow}(\mu)$ $|$ $\mu$ is open in $L_{S}$ and $f\in\mathscr F\}>$ and hence $\tau_{s}=<\mathscr F>$. Now we show that $\mathscr F$ is an $L$-topology on $X$. Let $\{f_{i}$ $|$ $i\in I\}\subseteq\mathscr F$. Define a map $g:(X,s)\rightarrow (S,u)^{I}$ as $g(x)(i)=f_{i}(x)$, for every $x\in X$ and for every $i\in I$. Let $p_{i}$ denote the $i$-th projection map from $S^{I}$ to $U(S,u)$. Then for every $i\in I$ and for every $x\in X$, $(p_{i}\circ g)(x)=p_{i}(g(x))=g(x)(i)=f_{i}(x)$ implying that $p_{i}\circ g=f_{i}$, for every $i\in I$. By condition $(2)$, $g$ is $\mathscr C$-morphism. For $x\in X$, $(sup\circ g)(x)=sup (g(x))=\displaystyle\bigvee_{i\in I}g(x)i=\displaystyle\bigvee_{i\in I}f_{i}(x)$ implying that $sup\circ g=\displaystyle\bigvee_{i\in I}f_{i}$. By using condition $(3)$, $\displaystyle\bigvee_{i\in I}f_{i}$ is $\mathscr C$-morphism and thus $\displaystyle\bigvee_{i\in I}f_{i}\in\mathscr F$. Similarly $\displaystyle\bigwedge_{i\in I}f_{i}\in\mathscr F$, for every finite $I$. Hence $\mathscr F$ is an $L$-topology on $X$. Thus $<\mathscr F>=\mathscr F$ and therefore $\tau_{s}=\mathscr F$. Let $f:(X,s)\rightarrow (Y,w)$ be $\mathscr C$-morphism. We show that $f:(X,\tau_{s})\rightarrow (Y,\tau_{w})$ is continuous. Let $\mu\in\tau_{w}$. Then $\mu:(Y,w)\rightarrow (S,u)$ is $\mathscr C$-morphism and hence $\mu\circ f:(X,s)\rightarrow (S,u)$ is also $\mathscr C$-morphism. As $f^{\leftarrow}(\mu)=\mu\circ f$, $f^{\leftarrow}(\mu)\in\tau_{s}$. Thus $f$ is continuous.\\
Let $(X,\tau)\in ob L\textbf{-Top}$ and $s_{\tau}$ be the initial lift of the family $\{\mu:X\rightarrow U(S,u)$ $|$ $\mu\in\tau\}$. Then $(X,s_{\tau})\in ob\mathscr C$ and $\{\mu:(X,s_{\tau})\rightarrow (S,u)$ $|$ $\mu\in\tau\}$ is an initial family.
Let $f:(X,\tau)\rightarrow(Y,\delta)$ be continuous. We show that $f:(X,s_{\tau})\rightarrow (Y,s_{\delta})$ is $\mathscr C$-morphism. For every $\nu\in\delta$, $f^{\leftarrow}(\nu)=\nu\circ f$ and hence $\nu\circ f\in\tau$. As $\nu\in\delta$, $\nu:(Y,s_{\delta})\rightarrow (S,u)$ is $\mathscr C$-morphism and $\nu\circ f\in\tau$, $f:(X,s_{\tau})\rightarrow (S,u)$ is $\mathscr C$-morphism. Since $\{\nu:(Y,s_{\delta})\rightarrow (S,u)$ $|$ $\nu\in\delta\}$ is an initial family, $f$ is $\mathscr C$-morphism.\\
Now we define two functor $F:\mathscr C\rightarrow L\textbf{-Top}$ as $F(X,s)=(X,\tau_{s})$ and $G:L\textbf{-Top}\rightarrow\mathscr C$ as $G(X,\tau)=(X,s_{\tau})$ (both functors have morphisms unchanged at set-theoretic level). It can be verified easily that $F$ and $G$ are concrete functors and they are inverses to each other. $\Box$
\section{Epireflective hull of $L_{S}$ in $L$-Top}
\qquad The category $\textbf{Top}_{0}$ of $T_{0}$-topological spaces is well-known to be the epireflective hull of usual two-point Sierpinski space $2_{S}$ in the category $\textbf{Top}$ of topological spaces (\cite{AHS}, page 263). Lowen and Srivastava \cite{LS} showed that the category $\textbf{FTop}_{0}$ of $T_{0}$-fuzzy topological spaces is the epireflective hull of fuzzy Sierpinski space $I_{S}$ in the category $\textbf{FTop}$ of fuzzy topological spaces. Analogously, in this section we show that $L\textbf{-Top}_{0}$ is the epireflective hull of $L_{S}$ in $L\textbf{-Top}$.

For showing that $L\textbf{-Top}_{0}$ is the epireflective hull of $L_{S}$ in $L\textbf{-Top}$, we first need to identify the epimorphisms and extremal subobjects in $L\textbf{-Top}$. Epimorphisms and extremal subobjects, along with their proofs in $L\textbf{-Top}$ are on familiar lines (as in $\textbf{Top}$).
\begin{pro}
Epimorphisms in $L\textbf{-Top}$ are precisely the surjective maps.
\end{pro}
\begin{pro}
Extremal monomorphisms in $L\textbf{-Top}$ are precisely the embeddings in $L\textbf{-Top}$.
\end{pro}

\begin{cor}
Extremal subobjects in $L\textbf{-Top}$ are precisely the subspaces.\\
\end{cor}

We Shall use the following two results from \cite{Mar} (Theorem $1$ and Theorem $2$) to obtaining the epireflective hull of $L_{S}$ in the categories $L\textbf{-Top}$ and $L\textbf{-Top}_{0}$.
\begin{thm} {\rm \cite{Mar}}
A subcategory $\mathscr B$ of a category $\mathscr A$ is epireflective in $\mathscr A$ iff it is closed under the formation of products and extremal subobjects in $\mathscr A$.
\end{thm}
\begin{thm} {\rm \cite{Mar}}
Let $\mathscr A$ be a category, $E$ be a class of $\mathscr A$-objects and $RE$ be the epireflective hull of $E$ in $\mathscr A$. Then $A\in ob RE$ iff $A$ is an extremal subobject of a product of objects of $E$.
\end{thm}

In the above two results, it is assumed that the category $\mathscr A$ has products and $\mathscr A$ is an epi-co-well-powered (epi, extremal mono)-category (cf. \cite{Mar}).\\

\begin{pro}
$L\textbf{-Top}$ is an epi-co-well-powered $($Epi, Extremal mono$)$-category.
\end{pro}
\textbf{Proof}: As pointed out earlier, $L\textbf{-Top}$ is a topological construct. Hence by using Theorem $21.16$ and Corollary $21.17$ of \cite{AHS}, the result follows. $\Box$
\begin{thm}
$L\textbf{-Top}_{0}$ is epireflective in $L\textbf{-Top}$.
\end{thm}
\textbf{Proof}: Let $(X,\tau)\in ob L\textbf{-Top}$. Define a relation $\sim$ on $X$ as follows: for every $x,y\in X$, $x\sim y$ iff $\mu(x)=\mu(y)$, for every $\mu\in\tau$. It is easy to see that $\sim$ is an equivalence relation on $X$. Consider the quotient map $q_{X}:X\rightarrow \tilde{X}$, where $\tilde X=X/ \sim $ and let $\tilde{\tau}$ be the quotient $L$-topology on $\tilde{X}$ induced by $(X,\tau)$ and $q_{X}$. Then $(\tilde X, \tilde \tau)\in ob L\textbf{-Top}_{0}$. It can be verified easily that the quotient map $q_X:(X,\tau)\rightarrow (\tilde X, \tilde \tau)$ is an epireflection of $(X,\tau)$ in $L\textbf{-Top}_{0}$. $\Box$\\

By using Theorem $5.1$, we get the following corollary.
\begin{cor}
$L\textbf{-Top}_{0}$ is closed under forming products and extremal subobjects in $L\textbf{-Top}$.
\end{cor}

\begin{pro}
$(X,\tau)\in ob L\textbf{-Top}_{0}$ iff it is homeomorphic to a subspace of product of copies of $L_{S}$.
\end{pro}
\textbf{Proof}: Let $(X,\tau)\in ob L\textbf{-Top}_{0}$. Define $e:(X,\tau)\rightarrow L_{S}^{\tau}$ as $e(x)(\mu)=\mu(x)$, for every $x\in X$ and for every $\mu\in\tau$. Let $p_{\mu}$ denote the $\mu$-th projection map from $L_{S}^{\tau}$ to $L_{S}$. Then for every $x\in X$, $(p_{\mu}\circ e)(x)=p_{\mu}(e(x))=e(x)(\mu)=\mu(x)$ implying that $p_{\mu}\circ e=\mu$. Thus $e$ is continuous. For distinct $x,y\in X$, there exists $\mu\in\tau$ such that $\mu(x)\neq \mu(y)$. Hence $e(x)\neq e(y)$, showing that $e$ is injective. Let $\mu\in\tau$. Then for every $x\in X$, $e^{\rightarrow}(\mu)e(x)=\vee\{\mu(y)$ $|$ $e(y)=e(x)\}=\mu(x)=e(x)(\mu)=p_{\mu}(e(x))$ implying that $e^{\rightarrow}(\mu)=p_{\mu}\wedge 1_{e(X)}$. As $p_{\mu}$ is open in $L_{S}^{\tau}$, $e^{\rightarrow}(\mu)$ is open in $e(X)$. Thus $(X,\tau)$ is homeomorphic to a subspace of product of copies of $L_{S}$. By using the Corollary $5.1$, Corollary $5.2$ and the fact that $L_{S}$ is $T_{0}$, the converse follows. $\Box$\\

In view of Theorem $5.2$, Corollary $5.1$, and Proposition $5.4$ above, we get the following result:
\begin{thm}
$L\textbf{-Top}_{0}$ is the epireflective hull of $L_{S}$ in $L\textbf{-Top}$.
\end{thm}
\section{The sober $L$-topological space}
\qquad Rodabaugh \cite{Roda} and H\"{o}hle \cite{Hole} introduced the notion of sober $L$-topological space in the category $L\textbf{-Top}$ as follows.\\

Let $A$ be a frame and let $pt_{L}A=\{p:A\rightarrow L$ $|$ $p$ is a frame map$\}$. Define $\phi_{L}:A\rightarrow L^{pt_{L}A}$ as $\phi_{L}(a)(p)=p(a)$, for every $a\in A$ and for every $p\in pt_{L}A$. Then $\phi_{L}(A)$ is an $L$-topology on $pt_{L}A$ (cf. \cite{Roda}).\\
For an $L$-topological space $(X,\tau)$, $\tau$ is a frame.
\begin{def1}{\rm \cite{Hole}}
An $L$-topological space $(X,\tau)$ is  said to be {\rm sober} if for every $p\in pt_{L}\tau$, there exists a unique $x\in X$ such that $p(\mu)=\mu(x)$, for every $\mu\in\tau$.
\end{def1}

Let $L\textbf{-Sob}$ denote the subcategory of $L\textbf{-Top}$, whose objects are sober $L$-topological spaces.

The following result is easy to verify.
\begin{pro}
The Sierpinski $L$-topological space $L_{S}$ is sober.
\end{pro}
\begin{pro}{\rm \cite{Roda}}
An $L$-topological space $(X,\tau)$ is sober iff $\eta_{X}:(X,\tau)\rightarrow (pt_{L}\tau, \phi_{L}(\tau))$ defined as $\eta_{X}(x)(\mu)=\mu(x)$, for every $x\in X$ and for every $\mu\in\tau$, is a homeomorphism.
\end{pro}
\begin{pro}{\rm \cite{Roda}}
For an $L$-topological space $(X,\tau)$, $(pt_{L}\tau, \phi_{L}(\tau))$ is sober.
\end{pro}

\begin{pro}{\rm \cite{Roda}}
An $L$-topological space $(X,\tau)$ is $T_{0}$ iff $\eta_{X}:(X,\tau)\rightarrow (pt_{L}\tau, \phi_{L}(\tau))$ defined as $\eta_{X}(x)(\mu)=\mu(x)$, for every $x\in X$ and for every $\mu\in\tau$, is injective.
\end{pro}
In view of Proposition $6.2$ and Proposition $6.4$, $L\textbf{-Sob}$ is a subcategory of $L\textbf{-Top}_{0}$.
\section{Epireflective hull of $L_{S}$ in $L$-Top$_{0}$}
\qquad It is known that $L\textbf{-Sob}$ is reflective in $L\textbf{-Top}$ (cf. \cite{{Roda}, {Hole}}). As $L\textbf{-Sob}$ is a subcategory of $L\textbf{-Top}_{0}$, $L\textbf{-Sob}$ is also reflective in $L\textbf{-Top}_{0}$. Here we show that $L\textbf{-Sob}$ is epireflective in $L\textbf{-Top}_{0}$. Moreover, we show that $L\textbf{-Sob}$ is, in fact, $\mathscr H$-firm epireflective in $L\textbf{-Top}_{0}$. For showing this, first we have to identify the epimorphisms in $L\textbf{-Top}_{0}$.\\

For $(X,\tau)\in ob L\textbf{-Top}$ and $M\subseteq X$, put $[M]=\cap\{Eq(f,g)$ $|$ $f,g:(X,\tau)\rightarrow (Y,\delta)$ are morphisms in $L\textbf{-Top}$ with $(Y,\delta)\in ob L\textbf{-Top}_{0}$ and $f|_{M}=g|_{M}\}$, where $Eq(f,g)=\{x\in X$ $|$ $f(x)=g(x)\}$. 
$M$ is said to be $[\quad]$-closed if $[M]=M$. It can be verified easily that $[[M]]=[M]$.

\begin{pro}
A morphism $f:(X,\tau)\rightarrow (Y,\delta)$ in $L\textbf{-Top}_{0}$ is an epimorphism iff for every $\nu_{1},\nu_{2}\in\delta$, $f^{\leftarrow}(\nu_{1})=f^{\leftarrow}(\nu_{2})$ implies that $\nu_{1}=\nu_{2}$.
\end{pro}
\textbf{Proof}: Let $f:(X,\tau)\rightarrow (Y,\delta)$ be an epimorphism in $L\textbf{-Top}_{0}$. Let $\nu_{1},\nu_{2}\in\delta$ such that
$f^{\leftarrow}(\nu_{1})=f^{\leftarrow}(\nu_{2})$. Then $\nu_{1}\circ f=\nu_{2}\circ f$, implying that $\nu_{1}=\nu_{2}$.

Conversely, let for every $\nu_{1},\nu_{2}\in\delta$, $f^{\leftarrow}(\nu_{1})=f^{\leftarrow}(\nu_{2})$ implies that $\nu_{1}=\nu_{2}$. Let $g,h:(Y,\delta)\rightarrow (Z,\Delta)$ be two distinct morphisms in $L$\textbf{-Top}$_{0}$. Then there exists $y\in Y$ such that $g(y)\neq h(y)$. As $(Z,\Delta)$ is $T_{0}$, there exists $\mu\in\Delta$ such that $\mu(g(y))\neq\mu(h(y))$. Thus $g^{\leftarrow}(\mu)\neq h^{\leftarrow}(\mu)$, implying that $f^{\leftarrow}(g^{\leftarrow}(\mu))\neq f^{\leftarrow}(h^{\leftarrow}(\mu))$, i.e., $\mu\circ g\circ f\neq \mu\circ h\circ f$. Hence $g\circ f\neq h\circ f$, showing that $f$ is an epimorphism. $\Box$\\

We shall use the following result, which is a special case of Theorem $1.11$ of \cite{Cast}.
\begin{pro}
A morphism $f:(X,\tau)\rightarrow (Y,\delta)$ in $L\textbf{-Top}_{0}$ is an epimorphism iff $[f(X)]=Y$.
\end{pro}


Call an embedding $e:(X,\tau)\rightarrow (Y,\delta)$ in $L\textbf{-Top}_{0}$, $[\quad]$-closed if $[e(X)]=e(X)$.
\begin{pro}
Extremal monomorphisms in $L\textbf{-Top}_{0}$ are precisely the $[\quad]$-closed embeddings.
\end{pro}
\textbf{Proof}: Let $f:(X,\tau)\rightarrow (Y,\delta)$ be an extremal monomorphism in $L\textbf{-Top}_{0}$. Let $f':(X,\tau)\rightarrow ([f(X)],\delta_{[f(X)]})$ be the corestriction of $f$ onto $[f(X)]$. Then $f'$ is an epimorphism. Now $f=i\circ f'$, where $i:([f(X)],\delta_{[f(X)]})\rightarrow (y,\delta)$ is the inclusion map, implying that $f'$ is an isomorphism. Thus $f$ is $[\quad]$-closed embedding.

Conversely, let $e:(X,\tau)\rightarrow (Y,\delta)$ be a $[\quad]$-closed embedding in $L\textbf{-Top}_{0}$. Then $[e(X)]=e(X)$ is an equalizer of some pair of morphism in $L\textbf{-Top}_{0}$(\cite{Cast}, Proposition $1.6$). But equalizers are necessarily extremal monomorphisms. So, $e$ is an extremal monomorphism. $\Box$
\begin{cor}
Extremal subobjects in $L\textbf{-Top}_{0}$ are precisely the $[\quad]$-closed subspaces.
\end{cor}
\begin{pro}
$L\textbf{-Top}_{0}$ is an epi-co-well-powered $($Epi, Extremal mono$)$-category.
\end{pro}
\textbf{Proof}: $L\textbf{-Top}_{0}$ is an epireflective subcategory of $L\textbf{-Top}$ (Theorem $5.3$). Also $L\textbf{-Top}$ is a topological category. By using Proposition $3.3$ of \cite{Giu} $L\textbf{-Top}_{0}$ is a well-powered and complete $($Epi, Extremal mono$)$-category. By using the Example $7.90$ $(2)$ of \cite{AHS}, $L\textbf{-Top}_{0}$ is an epi-co-well-powered and hence the result follows. $\Box$
\begin{thm}
$L\textbf{-Sob}$ is epireflective in $L\textbf{-Top}_{0}$.
\end{thm}
\textbf{Proof}: Let $(X,\tau)\in ob L\textbf{-Top}_{0}$. Then $\eta_{X}:(X,\tau)\rightarrow (pt_{L}\tau, \phi_{L}(\tau))$ defined as $\eta_{X}(x)(\mu)=\mu(x)$, for every $x\in X$ and for every $\mu\in\tau$, is the reflection of $(X,\tau)$ in $L\textbf{-Sob}$ (cf. \cite{Hole}, Corollary 2.1). Now, we show that $\eta_{X}$ is an epimorphism. Let $\mu_{1},\mu_{2}\in\tau$ and $\eta_{X}^{\leftarrow}(\phi_{L}(\mu_{1}))=\eta_{X}^{\leftarrow}(\phi_{L}(\mu_{2}))$. Then for every $x\in X$, $\eta_{X}^{\leftarrow}(\phi_{L}(\mu_{1}))(x)=\eta_{X}^{\leftarrow}(\phi_{L}(\mu_{2}))(x)$ implying that $\phi_{L}(\mu_{1})\eta_{X}(x)=\phi_{L}(\mu_{2})\eta_{X}(x)$. Then by definition of $\phi_{L}$, $\eta_{X}(x)(\mu_{1})=\eta_{X}(x)(\mu_{2})$ which implies that $\mu_{1}(x)=\mu_{2}(x)$, for every $x\in X$. Thus $\mu_{1}=\mu_{2}$ implying that $\phi_{L}(\mu_{1})=\phi_{L}(\mu_{2})$. Thus $\eta_{X}$ is an epimorphism. $\Box$\\

By using the Theorem $5.1$, we get the following corollary.
\begin{cor}
$L\textbf{-Sob}$ is closed under forming products and extremal subobjects in $L\textbf{-Top}_{0}$.\\
\end{cor}

For $(X,\tau)\in ob L\textbf{-Top}_{0}$, $\eta_{X}:(X,\tau)\rightarrow (pt_{L}\tau, \phi_{L}(\tau))$ defined as $\eta_{X}(x)(\mu)=\mu(x)$, for every $x\in X$ and for every $\mu\in\tau$, is an open map. As, for $\mu\in\tau$ and $p\in pt_{L}\tau$, $\eta_{X}^{\rightarrow}(\mu)(p)=\vee\{\mu(x)$ $|$ $\eta_{X}(x)=p\}=\vee\{\eta_{X}(x)(\mu)$ $|$ $\eta_{X}(x)=p\}=p(\mu)=\phi_{L}(\mu)(p)$, showing that $\eta_{X}^{\rightarrow}(\mu)\in \phi_{L}(\tau)$. Thus in view of Proposition $6.4$, $\eta_{X}$ is an embedding in $L\textbf{-Top}_{0}$.\\
\begin{thm}
$L\textbf{-Sob}$ is $\mathscr H$-firm epireflective in $L\textbf{-Top}_{0}$, where $\mathscr H$ denotes the class of all embeddings in $L\textbf{-Top}_{0}$.
\end{thm}
\textbf{Proof}: Let $(X,\tau)\in ob L\textbf{-Top}_{0}$. Then $\eta_{X}:(X,\tau)\rightarrow (pt_{L}\tau, \phi_{L}(\tau))$ is an epimorphic embedding in $L\textbf{-Top}_{0}$. Let $(Y,\delta)\in ob L\textbf{-Sob}$ and $f:(X,\tau)\rightarrow (Y,\delta)$ be an epimorphic embedding in $L\textbf{-Top}_{0}$. We have to find a unique isomorphism $f^{*}:(pt_{L}\tau, \phi_{L}(\tau))\rightarrow (Y,\delta)$ such that $f^{*}\circ \eta_{X}=f$. Let $p\in pt_{L}\tau$. Then $p:\tau\rightarrow L$ is a frame map. Define $p':\delta\rightarrow L$ as $p'(\nu)=p(\nu\circ f)$, for every $\nu\in\delta$. It is easy to see that $p'$ is a frame map. As $(Y,\delta)$ is sober, there exists a unique $y\in Y$ such that for every $\nu\in\delta$, $p'(\nu)=\nu(y)$ implying that $p(\nu\circ f)=\nu(y)$. Put $f^{*}(p)=y$. Let $\nu\in\delta$. Then for every $p\in pt_{L}\tau$, ${f^{*}}^{\leftarrow}(\nu)(p)=\nu(f^{*}(p))=p(\nu\circ f)=p(f^{\leftarrow}(\nu))=\phi_{L}(f^{\leftarrow}(\nu))(p)$, implying that ${f^{*}}^{\leftarrow}(\nu)=\phi_{L}(f^{\leftarrow}(\nu))$. As $\nu\in\delta$, $f^{\leftarrow}(\nu)\in\tau$ showing that ${f^{*}}^{\leftarrow}(\nu)\in\phi_{L}(\tau)$. Hence $f^{*}$ is continuous. For every $x\in X$, $(f^{*}\circ\eta_{X})(x)=f^{*}(\eta_{X}(x))$. Then for every $\nu\in\delta$, $\nu(f^{*}(\eta_{X}(x)))=\eta_{X}(x)(\nu\circ f)=(\nu\circ f)(x)=\nu(f(x))$. By the uniqueness of $y$ in the definition of $f^{*}$, $f^{*}(\eta_{X}(x))=f(x)$, for every $x\in X$. Hence $f^{*}\circ\eta_{X}=f$. Now, we show that $f^{*}$ is an isomorphism.\\
For every $\mu\in\tau$, $f^{\rightarrow}(\mu)\in \delta_{f(X)}$ implying that there exists $\mu_{f}\in\delta$ such that $f^{\rightarrow}(\mu)=\mu_{f}\wedge 1_{f(X)}$ and $f^{\leftarrow}(\mu_{f})=\mu$. As $f$ is an epimorphism, $\mu_{f}$ is unique such that $f^{\leftarrow}(\mu_{f})=\mu$. Define $g:(Y,\delta)\rightarrow (pt_{L}\tau,\phi_{L}(\tau))$ as $g(y)(\mu)=\mu_{f}(y)$, for every $y\in Y$ and for every $\mu\in\tau$. It can be verified easily that $g(y)$ is a frame map. Now we show that $g$ is continuous. Let $\mu\in\tau$. Then for every $y\in Y$, $g^{\leftarrow}(\phi_{L}(\mu))(y)=\phi_{L}(\mu)g(y)=g(y)(\mu)=\mu_{f}(y)$, showing that $g^{\leftarrow}(\phi_{L}(\mu))=\mu_{f}\in\delta$. Thus, $g$ is continuous.\\
For every $\nu\in\delta$, $f^{\leftarrow}(\nu)\in\tau$ implying that $f^{\leftarrow}((f^{\leftarrow}(\nu))_{f})=f^{\leftarrow}(\nu)$ and hence $(f^{\leftarrow}(\nu))_{f}=\nu$. For every $y\in Y$, $(f^{*}\circ g)(y)=f^{*}(g(y))$. Then for every $\nu\in\delta$, $\nu(f^{*}(g(y)))=g(y)(\nu\circ f)=g(y)(f^{\leftarrow}(\nu))=(f^{\leftarrow}(\nu))_{f}(y)=\nu(y)$. by the uniqueness of $y$ in the definition of $f^{*}$, $f^{*}(g(y))=y$, for every $y\in Y$. Hence $f^{*}\circ g=id_{Y}$.\\
For every $p\in pt_{L}\tau$ and $\mu\in\tau$, $(g\circ f^{*})(p)(\mu)=g(f^{*}(p))(\mu)=\mu_{f}(f^{*}(p))=p(\mu_{f}\circ f)=p(f^{\leftarrow}(\mu_{f}))=p(\mu)$ implying that $g\circ f^{*}=id_{pt_{L}\tau}$. $\Box$\\

It is well-known that the category $\textbf{Sob}$ of sober topological spaces is the epireflective hull of usual two-point Sierpinski space $2_{S}$ in the category $\textbf{Top}_{0}$ of $T_{0}$-topological spaces \cite{Nel}. Srivastava and Khastgir \cite{SK} showed that the category $\textbf{FSob}$ of fuzzy sober spaces is the epireflective hull of fuzzy Sierpinski space $I_{S}$ in the category $\textbf{FTop}_{0}$ of $T_{0}$-fuzzy topological spaces. Analogously, we show that $L\textbf{-Sob}$ is the epireflective hull of $L_{S}$ in $L\textbf{-Top}_{0}$.

\begin{pro}
For every $(X,\tau)\in ob L\textbf{-Top}_{0}$, $(pt_{L}\tau,\phi_{L}(\tau))$ is a $[\quad]$-closed subspace of $L_{S}^{\tau}$.
\end{pro}
\textbf{Proof}: 
Consider the epimorphic embedding $\eta_{X}:(X,\tau)\rightarrow (pt_{L}\tau,\phi_{L}(\tau))$ in $L\textbf{-Top}_{0}$. As $(X,\tau)\in obL\textbf{-Top}_{0}$, $(X,\tau)$ can be embedded in $L_{S}^{\tau}$ via the `evaluation map' $e:(X,\tau)\rightarrow L_{S}^{\tau}$ defined as $e(x)(\mu)=\mu(X)$, for every $x\in X$ and for every $\mu\in\tau$. Consider the subspace $[e(X)]$ of $L_{S}^{\tau}$. Let $e'$ be the corestriction of $e$ onto $e[X]$. Then $e'$ is an epimorphic embedding in $L\textbf{-Top}_{0}$. By using Proposition $6.1$, Corollary $7.1$ and Corollary $7.2$, $[e(X)]$ is sober. Then there exists an isomorphism from $(pt_{L}\tau,\phi_{L}(\tau))$ to $[e(X)]$ (by Theorem $7.2$). Hence $(pt_{L}\tau,\phi_{L}(\tau))$ is a $[\quad]$-closed subspace of $L_{S}^{\tau}$. $\Box$
\begin{pro}
$(X,\tau)\in ob L\textbf{-Sob}$ iff it is homeomorphic to a $[\quad]$-closed subspace $L_{S}^{\tau}$.
\end{pro}
\textbf{Proof}: Let $(X,\tau)\in ob L\textbf{-Sob}$. Then $\eta_{X}:(X,\tau)\rightarrow (pt_{L}\tau,\phi_{L}(\tau))$ is a homeomorphism. Hence $(X,\tau)$ is homeomorphic to a $[\quad]$-closed subspace $L_{S}^{\tau}$. By using Proposition $6.1$, Corollary $7.1$ and Corollary $7.2$, converse follows. $\Box$\\

In view of Theorem $5.2$, Corollary $7.1$ and Proposition $7.6$, we get the following result:
\begin{thm}
$L\textbf{-Sob}$ is the epireflective hull of $L_{S}$ in $L\textbf{-Top}_{0}$.\\\\\\
\end{thm}
\textbf{Acknowledgement}: The first author would like to thank the University Grants Commission, New Delhi, India, for financial support through its Senior Research Fellowship.

\end{document}